# WIDTH AND MODE OF THE PROFILE FOR SOME RANDOM TREES OF LOGARITHMIC HEIGHT

By Luc Devroye[1] and Hsien-Kuei Hwang[1,2]

*McGill University and Academia Sinica*

We propose a new, direct, correlation-free approach based on central moments of profiles to the asymptotics of width (size of the most abundant level) in some random trees of logarithmic height. The approach is simple but gives precise estimates for expected width, central moments of the width and almost sure convergence. It is widely applicable to random trees of logarithmic height, including recursive trees, binary search trees, quad trees, plane-oriented ordered trees and other varieties of increasing trees.

**1. Introduction.** Most random trees in the discrete probability literature have height either of order $\sqrt{n}$ or of order $\log n$ ($n$ being the tree size); see [1]. For simplicity, we call these trees *square-root trees* and *log trees*, respectively. Profiles (number of nodes at each level of the tree) of random square-root trees have a rich connection to diverse structures in combinatorics and in probability, and have been extensively studied. In contrast, profiles of random log trees, arising mostly from data structures and computer algorithms, were less addressed and only quite recently were their limit distributions, drastically different from those of square-root trees, better understood; see [3, 12, 13, 21, 27] and the references therein.

We study in this paper the asymptotics of width, which is defined to be the size of the most abundant level, and its close connection to the profile. There are many results on first-order asymptotics of profiles for standard log trees, such as binary search trees, random recursive trees, $m$-ary search trees and quad trees. In some cases, quite accurate asymptotic expressions are known for the expected profile. There is already a paucity of results with regard to

Received August 2005; revised October 2005.
[1]Supported in part by the Humboldt Foundation.
[2]Supported in part by a grant from National Science Council of the Republic of China.

higher moments of the profile, let alone limit laws for properly centered and scaled profiles. The literature on this subject was surveyed by Drmota and Hwang [12, 13] and Fuchs, Hwang and Neininger [21], but some key historical references are repeated here. In fact, the last paper describes the complete asymptotics for the profile of random recursive trees and random binary search trees, two important representatives. The present paper extends the results to a universal class of log trees called *width-regular*. It adds a host of new results on the width of these trees, such as its exact location up to $O(1)$ as well as tight estimates on the central moments of the width and some strong laws of large numbers. Equally important is the fact that these results are obtained from simple moment estimates. For example, the correlation between the profiles at different levels is not needed at all in the study of the width.

*Recursive trees.* A prototypical log tree is the recursive tree, which has been introduced in diverse fields due to its simple construction. We will present our methods of proof for recursive trees and then indicate the required elements needed for other random trees.

*Combinatorially*, recursive trees are rooted, labeled, nonplane trees such that the labels along any path down from any node form an increasing sequence. By *random recursive trees*, we assume that all recursive trees of $n$ nodes are equally likely. *Probabilistically*, they can be constructed by successively adding nodes as follows. Start from a single root node with label 1. Then at the $i$th stage, the new node with label $i$ chooses any of the previous $i-1$ nodes uniformly at random [each with probability $1/(i-1)$] and is then attached to that node. This construction implies that there are $(n-1)!$ recursive trees of size $n$. The first paper on such tree models we could find is Tapia and Myers (under the name of concave node-weighted trees); see [45] and [13, 21, 43] for more references on the literature of recursive trees and their uses in other fields.

Note that the term "recursive trees" (first used in [35] and [37]) is less specific and has also been used in different contexts for different objects. For example, they are used in recursion computation theory to represent a computable set of strings with branching structure and in compilers to record the history of recursive procedures. They also appeared in classification trees, dynamic systems and data base languages with a different meaning.

*Profile.* Let $Y_{n,k}$ denote the number of nodes at distance $k$ from the root in random recursive trees of $n$ nodes (the root being at level zero). Such a profile is very informative and closely related to many other shape parameters, although it does not uniquely characterize the tree.

The combinatorial sister of the random recursive tree is the random binary search tree. Early work by Lynch [30] and Knuth [28] (see also [42]) showed that the expected profile of random binary search trees is related to Stirling numbers of the first kind and it peaks at $k$ about $2 \log n$. It peaks at $\log n + O(1)$ in random recursive trees. In contrast, the random variable $Y_{n,k}$ has received less attention until recently.

The profiles of random binary search trees and random recursive trees exhibit many interesting phenomena such as (i) bimodality of the variance, (ii) different ranges for convergence in distribution and for convergence of all moments of the normalized profile $Y_{n,k}/\mathbb{E}\{Y_{n,k}\}$, (iii) no convergence to a fixed limit law at the peak levels and (iv) sharp sign changes for the correlation coefficients of two level sizes; see [13, 21] for more information. See also [27], where the limit distribution of the profile of random binary search trees was first studied, and [3], where the width was addressed. However, their approach, which is based on martingale arguments, is very different from the moment approach used in this paper.

For simplicity, write throughout this paper $L_n := \max\{\log n, 1\}$. The expected profile $\mu_{n,k} := \mathbb{E}\{Y_{n,k}\}$, which gives the first picture of the general silhouette of random recursive trees, is known to be enumerated by the signless Stirling numbers of the first kind (see [43] and [21])

$$\sum_k \mu_{n,k} u^k = \prod_{1 \le j < n} \left(1 + \frac{u}{j}\right) = \binom{n+u-1}{n-1}.$$

From this, it follows by the saddle point method that

$$(1.1) \qquad \mu_{n,k} = \frac{n}{\sqrt{2\pi L_n}} e^{-\Delta^2/(2L_n) + O(|\Delta|^3/L_n^2)} \left(1 + O\left(\frac{1+|\Delta|}{L_n}\right)\right)$$

uniformly for $k = L_n + O(L_n^{2/3})$, where, *here and throughout this paper*, $\Delta := k - L_n$. The asymptotic approximation (1.1) is crucial for our analysis. It means that most nodes in a random recursive tree are located at the levels with $k \sim L_n$. In particular, we have

$$(1.2) \qquad \max_k \mu_{n,k} = \frac{n}{\sqrt{2\pi L_n}}(1 + O(L_n^{-1}));$$

see [46] or [25] for more precise expansions for $\mu_{n,k}$.

*Expected width.* We define the width of random recursive trees to be $W_n := \max_k Y_{n,k}$. The approximation (1.1) can be interpreted as a local limit theorem for the *depth*, which is the distance to the root of a randomly chosen node in recursive trees (each with the same probability). Thus it is intuitively clear that the width will be roughly close to $n/\sqrt{2\pi L_n}$. Our first result gives a more precise description of this.

THEOREM 1.1. *The expected width satisfies*

$$(1.3) \qquad \mathbb{E}\{W_n\} = \frac{n}{\sqrt{2\pi L_n}}(1 + \Theta(L_n^{-1})).$$

This theorem improves upon the error term $O(L_n^{-1/4} \log L_n)$ given in [13], where the proof depends on estimates for correlations of two level sizes and a tightness argument for the process. The approximation (1.3) also says, when compared with (1.1), that

$$\mathbb{E}\{W_n\} = \mu_{n,L_n+O(1)}(1 + O(L_n^{-1})).$$

In particular, by (1.2),

$$\mathbb{E}\left\{\max_k Y_{n,k}\right\} = \max_k \mathbb{E}\{Y_{n,k}\}(1 + O(L_n^{-1})).$$

Note that the index $\hat{k}$ reaching the maximum of $\mu_{n,k}$ satisfies

$$(1.4) \qquad\qquad \hat{k} = \lfloor L_n - 1 + \gamma + O(L_n^{-1}) \rfloor;$$

see [22] or pages 140–141 of [24]. Erdős [14] showed that $\hat{k}$ is unique.

*An estimate for absolute central moments.* We see from (1.3) that the expected width is asymptotically of the same order as the expected level sizes at $k = L_n + O(1)$. We show that not only their expected values are of the same order, but also all higher absolute central moments are asymptotically close.

THEOREM 1.2. *For any $s \geq 0$,*

$$(1.5) \qquad\qquad \mathbb{E}\{|W_n - \mathbb{E}\{W_n\}|^s\} = O(n^s L_n^{-3s/2}).$$

From [21], we have

$$(1.6) \qquad \mathbb{E}\{(Y_{n,k} - \mu_{n,k})^m\} = O(|\Delta|^m L_n^{-m} \mu_{n,k}^m) \qquad (k = L_n + o(L_n)).$$

By Lyapounov's inequality (see page 174 of [29]), we obtain, for any $s \geq 0$,

$$(1.7) \qquad \mathbb{E}\{|Y_{n,k} - \mu_{n,k}|^s\} = O(|\Delta|^s L_n^{-s} \mu_{n,k}^s) \qquad (k = L_n + o(L_n)).$$

In particular, this implies, by (1.1), that

$$\mathbb{E}\{|W_n - \mathbb{E}\{W_n\}|^s\} = O(\mathbb{E}\{|Y_{n,k} - \mu_{n,k}|^s\}) \qquad (s \geq 0)$$

for $k = L_n + O(1)$.

*Almost sure convergence.* In sequential growth models of random trees, such as the incremental model for random recursive trees described above, it makes sense to study almost sure convergence properties. One can also define sequential growth in random binary search trees and most other log trees discussed in this work. This is because most definitions derive from data storage applications, where the sequential insertion of new data in trees is a natural way to grow them. For random recursive trees, we use (1.5) to show that

$$(1.8) \qquad \frac{W_n}{\mathbb{E}\{W_n\}} \longrightarrow 1 \qquad \text{almost surely.}$$

This result was proved in [13] by martingale arguments and complex analysis, following [3]. Our proof relies on (1.5) with $s = 2 + \varepsilon$ and the usual Borel–Cantelli argument. It is conceptually simpler and also applies to random trees for which no martingale structure is available.

*Level reaching the width.* Let $k^*$ denote one of the levels such that $Y_{n,k^*} = W_n$. We show that $k^*$ takes most likely the values $L_n + O(1)$.

THEOREM 1.3. *For every $B > 0$, there exists $T_0 > 1$ such that*

$$\mathbb{P}(|k^* - L_n| \geq T) = O(T^{-B})$$

*for $T > T_0$.*

Thus width will with very small probability lie outside the range $L_n + O(1)$.

*Approaches used.* The most notable feature of our method of proof is that with the two crucial estimates (1.1) and (1.6) at hand, only basic probability tools such as Markov and Chebyshev inequalities and the Borel–Cantelli lemma are used. However, asymptotic tools for proving the two estimates for general random trees may differ from one case to another. In most cases we considered, the estimate (1.1) is proved by a combination of diverse analytic tools such as differential equations, singularity analysis (see [19]) and the saddle point method. The remaining analysis required for higher central moments of the profile is then mostly elementary, since this corresponds roughly to the large "toll-functions" cases for the underlying bivariate recurrences; see [21]. Although tools for handling the width of random square-root trees are very different from those for random log trees considered in this paper, the profiles provide in both classes of trees an accessible route to the asymptotics of the width; see [11] and the references therein.

*Generality of the phenomena.* The treatment for random recursive trees can be extended to *width-regular trees*, a large class of log trees defined below. This class includes familiar trees such as random $m$-ary search trees, quad trees, grid trees and increasing trees. To check whether a tree is width-regular is done case by case, unfortunately, because we need uniform estimates on the expected profile. A universal asymptotic tool is still lacking to extend the results, for example, to all random split trees [8].

*Organization of the paper.* For self-containedness and to pave the way for general random trees, we give a sketch of the proof for (1.1) and (1.6) in the next section. We then prove the theorems in Section 3. Extension of the same arguments to other log trees is given in Sections 4–7.

NOTATION. Throughout this paper, the generic symbol $\varepsilon > 0$ always represents a sufficiently small constant whose value may differ from one occurrence to another. Also $L_n := \max\{\log n, 1\}$ and $\Delta := k - L_n$.

**2. Estimates for the profile moments.** We briefly sketch the main ideas that lead to the estimates (1.1) and (1.6); see [21] for details and more precise estimates than (1.6).

*Recurrence of $Y_{n,k}$.* By construction, the profile of random recursive trees satisfies the recurrence

$$Y_{n,k} \overset{d}{=} \sum_{1 \le s < n} \frac{1}{s!} \underbrace{\sum_{\substack{j_1 + \cdots + j_s = n-1 \\ j_1, \ldots, j_s \ge 0}} \binom{n-1}{j_1, \ldots, j_s} \frac{(j_1-1)! \cdots (j_s-1)!}{(n-1)!}}_{\mathbb{P}\left(\substack{\text{the root degree equals } s \text{ and} \\ \text{the } s \text{ subtrees have sizes } j_1, \ldots, j_s}\right)}$$

$$\times (Y_{j_1,k-1}^{(1)} + \cdots + Y_{j_s,k-1}^{(s)})$$

for $n \ge 2$ and $k \ge 1$ with $Y_{1,0} = 1$, where the $Y_{n,k}^{(i)}$'s are independent copies of $Y_{n,k}$. From this we deduce, by conditioning on the size of the first subtree, that

$$(2.1) \qquad Y_{n,k} \overset{d}{=} Y_{I_n,k-1} + Y_{n-I_n,k}^* \qquad (n \ge 2; k \ge 1),$$

with $Y_{1,0} = 1$, where the $Y_{n,k}^*$'s are independent copies of $Y_{n,k}$ and independent of $I_n$, which is uniformly distributed in $\{1, \ldots, n-1\}$.

*The expected profile and the expansion* (1.1). From (2.1), we derive, by taking expectation and by solving the resulting recurrence, the relation

$$\sum_k \mu_{n,k} u^k = \binom{n+u-1}{n-1} \qquad (u \in \mathbb{C});$$

see [10, 36, 43]. Then by singularity analysis (see [19]),

$$(2.2) \qquad \sum_k \mu_{n,k} u^k = \frac{n^u}{\Gamma(1+u)}(1 + O(|u|(1+|u|)n^{-1})),$$

where the $O$-term holds uniformly for $|u| \le C$ for any $C > 0$. Note that the Stirling formula with complex parameter for the Gamma function does not give the required uniformity in $u$.

The uniform approximation (1.1) is then obtained by Cauchy's integral formula using (2.2) and the saddle point method.

*A uniform estimate for $\mu_{n,k}$.* A very useful uniform estimate for $\mu_{n,k}$ is given by

$$(2.3) \qquad \mu_{n,k} = O(L_n^{-1/2} r^{-k} n^r) \qquad (0 < r = O(1))$$

uniformly for all $0 \le k \le n$. This is easily obtained by Cauchy's integral formula and (2.2) since

$$\mu_{n,k} = O\left( r^{-k} n^r \int_{-\pi}^{\pi} n^{-r(1-\cos t)} \, dt \right),$$

which gives (2.3). *Throughout this paper, $r$ is always taken to be $r = 1 + o(1)$ unless otherwise specified.*

Although one can prove that

$$(2.4) \qquad \mu_{n,k} = O\left( \frac{L_n^k}{k!} \right)$$

uniformly for $0 \le k \le n$, the reason for using the estimate (2.3) instead of (2.4) is that, for general random search trees, it is much harder to derive the Poisson type estimate (2.4) for all $k$.

*Recurrence of higher central moments.* Let $P_{n,k}^{(m)} = \mathbb{E}\{(Y_{n,k} - \mu_{n,k})^m\}$. Then $P_{n,k}^{(m)}$ satisfies, by (2.1), the recurrence

$$P_{n,k}^{(m)} = \frac{1}{n-1} \sum_{1 \le j < n} (P_{j,k-1}^{(m)} + P_{n-j,k}^{(m)}) + Q_{n,k}^{(m)},$$

with $P_{n,0}^{(m)} = 0$ for $n, m \ge 1$, where

$$Q_{n,k}^{(m)} := \sum_{(a,b,c) \in \mathcal{I}_m} \binom{m}{a, b, c} \frac{1}{n-1} \sum_{1 \le j < n} P_{j,k-1}^{(a)} P_{n-j,k}^{(b)} \nabla_{n,k}^c(j) \qquad (m \ge 2),$$

with $\nabla_{n,k}(j) := \mu_{j,k-1} + \mu_{n-j,k} - \mu_{n,k}$ and

$$\mathcal{I}_m := \{(a,b,c) \in \mathbb{Z}^3 : a+b+c = m, 0 \le a, b < m, 0 \le c \le m\}.$$

We prove (1.6) in two stages. A uniform estimate for $\nabla_{n,k}(j)$ for $1 \le j, k < n$ is first derived, which then implies by induction a uniform bound for $P_{n,k}^{(m)}$ for $1 \le k < n$. This bound is, however, not tight when $\Delta = o(\sqrt{L_n})$. Then we refine the estimate for $\nabla_{n,k}(j)$ when $\Delta = O(\sqrt{L_n})$, which then leads to (1.6) by another induction.

*First estimate for $P_{n,k}^{(m)}$.* By (2.2), we have the integral representation

$$(2.5) \qquad \begin{aligned} \nabla_{n,k}(j) = \frac{1}{2\pi i} &\oint_{|u|=r} \frac{u^{-k-1} n^u}{\Gamma(1+u)} \varphi(u; j/n) \\ &\times (1 + O(j^{-1} + (n-j)^{-1})) \, du, \end{aligned}$$

where $\varphi(u; x) := ux^u + (1 - x)^u - 1$. Since $\varphi(1; x) = 0$, we have

$$\frac{\varphi(u; x)}{\Gamma(1 + u)} = O(|u - 1|)$$

uniformly for $x \in [0, 1]$. Substituting this estimate into (2.5), we obtain

$$
\begin{aligned}
(2.6) \qquad \nabla_{n,k}(j) &= O\left(r^{-k} n^r \int_{-\pi}^{\pi} |r e^{i\theta} - 1|^{n - r(1 - \cos\theta)} \, d\theta\right) \\
&= O((|r - 1| + L_n^{-1/2}) L_n^{-1/2} r^{-k} n^r)
\end{aligned}
$$

uniformly for $1 \le j, k < n$, where $r = 1 + o(1)$. This bound is not tight for all $k$, but is sufficient for most of our purposes. In particular, since $r$ is not specially chosen to minimize the error term, (2.6) is not optimal when $|r - 1| = o(L_n^{-1/2})$, which is the case when we choose $r = k/L_n$ and $|k - L_n| = o(\sqrt{L_n})$.

We now prove by induction that

$$(2.7) \qquad P_{n,k}^{(m)} = O((|r - 1|^m + L_n^{-m/2}) L_n^{-m/2} r^{-km} n^{mr}) \qquad (m \ge 0)$$

uniformly for $1 \le k < n$.

Obviously, (2.7) holds for $m = 0, 1$. Assume $m \ge 2$. To estimate $Q_{n,k}^{(m)}$, we split the sum into two parts,

$$Q_{n,k}^{(m)} = \sum_{(a,b,c) \in \mathcal{I}_m} \binom{m}{a, b, c} \frac{1}{n - 1} \left(\sum_{j \in \mathcal{J}_m} + \sum_{j \in \mathcal{J}_m'}\right) P_{j,k-1}^{(a)} P_{n-j,k}^{(b)} \nabla_{n,k}^c(j),$$

where $\mathcal{J}_m := \{j : n/L_n^m \le j \le n - n/L_n^m\}$ and $\mathcal{J}_m' := \{1, \ldots, n - 1\} \setminus \mathcal{J}_m$. Then by induction and (2.6), the terms in $Q_{n,k}^{(m)}$ with $j \in \mathcal{J}_m'$ are bounded above by

$$
\begin{aligned}
O\Bigg(r^{-mk} n^{-1} \sum_{(a,b,c) \in \mathcal{I}_m} \Bigg( &L_n^{-(b+c)/2} n^{(b+c)r} \sum_{j < n/L_n^m} L_j^{-a/2} j^{ar} \\
&+ L_n^{-(a+c)/2} n^{(a+c)r} \sum_{j < n/L_n^m} L_j^{-b/2} j^{br} \Bigg)\Bigg) \\
= O(L_n^{-3m/2} r^{-mk} n^{mr})
\end{aligned}
$$

uniformly for $1 \le k < n$.

On the other hand, when $j \in \mathcal{J}_m$, we have $L_j \sim L_{n-j} \sim L_n$; thus by induction and the two estimates (2.6) and (2.7),

$$
\begin{aligned}
\sum_{(a,b,c) \in \mathcal{I}_m} \binom{m}{a, b, c} \frac{1}{n - 1} \sum_{j \in \mathcal{J}_m} P_{j,k-1}^{(a)} P_{n-j,k}^{(b)} \nabla_{n,k}^c(j) \\
= O((|r - 1|^m + L_n^{-m/2}) L_n^{-m/2} r^{-km} n^{mr}),
\end{aligned}
$$

it follows that

$$(2.8) \qquad Q_{n,k}^{(m)} = O((|r-1|^m + L_n^{-m/2})L_n^{-m/2}r^{-km}n^{mr})$$

uniformly for $1 \le k < n$.

From [21], we have the closed-form expression

$$(2.9) \qquad P_{n,k}^{(m)} = Q_{n,k}^{(m)} + \sum_{1 \le j < n} \sum_{0 \le \ell \le k} \frac{Q_{j,k-\ell}^{(m)}}{j}[u^\ell](u+1) \prod_{j < h < n} \left(1 + \frac{u}{h}\right),$$

where $[u^\ell]F(u)$ denotes the coefficient of $u^\ell$ in the Taylor expansion of $F$. Substituting the estimate (2.8), we obtain

$$P_{n,k}^{(m)} = O\Bigg( Q_{n,k}^{(m)} + r^{-km} \sum_{1 \le j < n} (|r-1|^m + L_j^{-m/2})L_j^{-m/2}j^{mr-1} $$
$$\times \sum_{0 \le \ell \le k} r^{m\ell}[u^\ell](u+1) \prod_{j < h < n} \left(1 + \frac{u}{h}\right) \Bigg).$$

Now

$$\sum_{0 \le \ell \le k} r^{m\ell}[u^\ell](u+1) \prod_{j < h < n} \left(1 + \frac{u}{h}\right) \le (1 + r^m) \prod_{j < h < n} \left(1 + \frac{r^m}{h}\right)$$
$$= O\left(\left(\frac{n}{j}\right)^{r^m}\right).$$

Thus (2.7) follows.

When $k \sim L_n$, we take $r = k/L_n$ in (2.7), giving

$$P_{n,k}^{(m)} = O((|\Delta|^m + L_n^{m/2})L_n^{-m}\mu_{n,k}^m),$$

which proves (1.6) when $\sqrt{L_n} \le |\Delta| = o(L_n)$.

*Proof of (1.6) when $\Delta = O(\sqrt{L_n})$.* We now refine the above procedure and prove (1.6) when $\Delta = O(\sqrt{L_n})$, which has the form

$$(2.10) \qquad P_{n,k}^{(m)} = O(|\Delta|^m L_n^{-3m/2} n^m) \qquad (m \ge 0).$$

By applying the expansion

$$\varphi(u;x) = \varphi_u'(1;x)(u-1) + O(|u-1|^2) \qquad (x \in [0,1])$$

and the usual saddle point method to (2.5), we deduce that

$$(2.11) \qquad \nabla_{n,k}(j) = O(|\Delta|L_n^{-3/2}n)$$

uniformly for $\Delta = O(\sqrt{L_n})$ and $1 \le j < n$. Note that this estimate also follows from (1.1).

We apply the same inductive procedure used to prove (2.8). By applying (2.7) to terms with $j \in \mathcal{J}_m'$ and (2.10) to terms with $j \in \mathcal{J}_m$ [starting from (2.11)], we have

$$Q_{n,k}^{(m)} = O(|\Delta|^m L_n^{-3m/2} n^m) \qquad (m \geq 2)$$

uniformly for $\Delta = O(\sqrt{L_n})$. This estimate and (2.8) gives, by (2.9), (2.10) and a similar decomposition of the sums involved,

$$P_{n,k}^{(m)} = O\Bigg( Q_{n,k}^{(m)} + \sum_{j \in \mathcal{J}_m} \sum_{0 \leq \ell = o(L_n)} |k - \ell - L_j|^m L_j^{-3m/2} j^{m-1} [u^\ell]$$

$$\times (u+1) \prod_{j < h < n} \left( 1 + \frac{u}{h} \right) \Bigg)$$

$$= O(|\Delta|^m L_n^{-3m/2} n^m).$$

This completes the proof of (2.10).

Such a two-stage proof of (1.6) is completely general when we have an integral representation for $\nabla_{n,k}(j)$ of the form (2.5) and a closed-form similar to (2.9). We will sketch means to handle the cases when no closed-form solution like (2.9) is available.

## 3. Asymptotics of the moments of the width.

We first prove Theorem 1.1; then we extend the proof for (1.5) and finally prove Theorem 1.3.

### 3.1. *Expected width.*

*Lower bound for the expected width.* The lower bound follows easily from the inequality

$$\mathbb{E}\{W_n\} \geq M_n,$$

where

$$M_n := \max_k \mathbb{E}\{Y_{n,k}\} = \frac{n}{\sqrt{2\pi L_n}}(1 + O(L_n^{-1}));$$

see (1.2).

*An inequality for the upper bound.* For the upper bound, we use the inequality

$$\mathbb{E}\{W_n\} \leq M_n + \sum_{|\Delta| \leq K} \mathbb{E}\{(Y_{n,k} - M_n)_+\} + \sum_{|\Delta| > K} \mu_{n,k}$$

$$=: w_n^{(1)} + w_n^{(2)} + w_n^{(3)},$$

where $K := L_n^{2/3}$.

*The sum $w_n^{(3)}$.* The last sum is easily estimated, since by (2.3),

$$w_n^{(3)} = O\left(L_n^{-1/2} n^r \left(\sum_{0 \le k \le L_n - K} + \sum_{k \ge L_n + K}\right) r^{-k}\right).$$

Taking $r = 1 - L_n^{-1/3}$, we see that

$$\sum_{0 \le k \le L_n - K} \mu_{n,k} = O\left(L_n^{-1/2} n^r \frac{r^{-L_n + K}}{1 - r}\right)$$

$$= O(L_n^{-1/2} n^{1 - L_n^{-1/3}} L_n^{1/3} (1 - L_n^{-1/3})^{-L_n + L_n^{2/3}})$$

$$= O(n L_n^{-1/6} e^{-L_n^{1/3}/2}).$$

The same upper bound holds for $\sum_{k \ge L_n + K} \mu_{n,k}$ by taking $r = 1 + L_n^{-1/3}$.

*An estimate for the second sum $w_n^{(2)}$.* We use the inequalities

$$\mathbb{E}\{(Y_{n,k} - M_n)_+\} \le \mathbb{E}\{(Y_{n,k} - \mu_{n,k}) \mathbf{1}_{(Y_{n,k} > M_n)}\}$$

$$\le \frac{\mathbb{E}\{(Y_{n,k} - \mu_{n,k})^2\}}{M_n - \mu_{n,k}}$$

for those $k$'s for which $M_n > \mu_{n,k}$. By (1.1),

$$
\begin{aligned}
M_n - \mu_{n,k} &= \frac{n}{\sqrt{2\pi L_n}} (1 - e^{-\Delta^2/(2L_n) + O(|\Delta|^3/L_n^2)})(1 + o(1)) \\
&= \frac{n}{\sqrt{2\pi L_n}} (1 - e^{-\Delta^2/(2L_n)} \\
&\quad + O(e^{-\Delta^2/(2L_n)} |\Delta|^3 L_n^{-2}))(1 + o(1)) \\
&\ge \frac{n}{\sqrt{2\pi L_n}} (1 - e^{-\Delta^2/(3L_n)})(1 + o(1))
\end{aligned}
\tag{3.1}
$$

uniformly for $1 \le |\Delta| \le K$. On the other hand, the variance is bounded above by

$$\mathbb{V}\{Y_{n,k}\} = O(\Delta^2 L_n^{-2} \mu_{n,k}^2) = O(\Delta^2 L_n^{-3} n^2 e^{-\Delta^2/L_n}) \tag{3.2}$$

uniformly for $1 \le |\Delta| \le K$. It follows from these estimates that

$$
\begin{aligned}
w_n^{(2)} &\le \sqrt{\mathbb{V}\{Y_{n,\lfloor L_n \rfloor}\}} + \sum_{1 \le |\Delta| \le K} \frac{\mathbb{V}\{Y_{n,k}\}}{M_n - \mu_{n,k}} \\
&= O(n L_n^{-3/2}) + O\left(L_n^{-5/2} n \int_1^\infty \frac{x^2 e^{-x^2/L_n}}{1 - e^{-x^2/(3L_n)}} \, dx\right) \\
&= O(n L_n^{-1}).
\end{aligned}
$$

Collecting all estimates, we get a weaker error term than (1.3),

$$(3.3) \qquad \mathbb{E}\{W_n\} = \frac{n}{\sqrt{2\pi L_n}}(1 + O(L_n^{-1/2})),$$

but we only used estimates for $\mathbb{E}\{Y_{n,k}\}$ and $\mathbb{V}\{Y_{n,k}\}$.

*Improving the error term by the fourth central moment of $Y_{n,k}$.* We can improve the error term in (3.3) by using the estimate for the fourth central moment of $Y_{n,k}$; see (1.6). Taking $m = 4$ in (1.6) and repeating the same analysis as above,

$$\begin{aligned}
w_n^{(2)} &\le \sqrt{\mathbb{V}\{Y_{n,\lfloor L_n \rfloor}\}} + \sum_{1 \le |\Delta| \le K} \frac{\mathbb{E}\{(Y_{n,k} - \mu_{n,k})^4\}}{(M_n - \mu_{n,k})^3} \\
&= O(nL_n^{-3/2}) + O\left(nL_n^{-9/2} \int_1^\infty \frac{x^4 e^{-2x^2/L_n}}{(1 - e^{-x^2/(3L_n)})^3}\, dx\right) \\
&= O(nL_n^{-3/2}) + O\left(nL_n^{-2} \int_{1/L_n}^\infty \frac{v^{3/2} e^{-2v}}{(1 - e^{-v/3})^3}\, dv\right) \\
&= O(nL_n^{-3/2}) + O\left(nL_n^{-2} \int_{1/L_n}^\infty v^{-3/2}\, dv\right) \\
&= O(nL_n^{-3/2}).
\end{aligned}$$

This proves (1.3).

3.2. *Higher absolute central moments of $W_n$.* We prove only an upper bound for $s = 2$, namely for the variance of $W_n$; other values of $s$ follow by the same argument and Lyapounov's inequality.

*An upper bound for the variance of the width.* We show, by using central moments of $Y_{n,k}$ of order 6, that

$$(3.4) \qquad \mathbb{V}\{W_n\} = O(n^2 L_n^{-3}),$$

which proves (1.5) with $s = 2$.

The proof extends that for $\mathbb{E}\{W_n\}$. Define $k_0 = \lfloor L_n \rfloor$. We start from

$$\begin{aligned}
\mathbb{E}\{(W_n - \mathbb{E}\{W_n\})^2\} &= \mathbb{E}\{(W_n - \mu_{n,k_0} + \mu_{n,k_0} - \mathbb{E}\{W_n\})^2\} \\
&\le 2\mathbb{E}\{(W_n - \mu_{n,k_0})^2\} + 2\mathbb{E}\{(\mu_{n,k_0} - \mathbb{E}\{W_n\})^2\}.
\end{aligned}$$

By (1.3),

$$\mathbb{E}\{(\mu_{n,k_0} - \mathbb{E}\{W_n\})^2\} = O(n^2 L_n^{-3})$$

and, similarly to the analysis for $\mathbb{E}\{W_n\}$,

$$
\begin{aligned}
\mathbb{E}\{(W_n - \mu_{n,k_0})^2\} &\leq \mathbb{E}\left\{\sum_{|\Delta| \geq 0} (Y_{n,k_0+\Delta} - \mu_{n,k_0})_+^2 \cdot \mathbf{1}_{(Y_{n,k_0+\Delta} > \mu_{n,k_0})}\right\} \\
&\leq \mathbb{V}\{Y_{n,k_0}\} + \sum_{1 \leq |\Delta| \leq K} \mathbb{E}\{(Y_{n,k_0+\Delta} - \mu_{n,k_0})_+^2\} \\
&\quad + \sum_{|\Delta| \geq K} \mathbb{V}\{Y_{n,k_0+\Delta}\} \\
&=: v_n^{(1)} + v_n^{(2)} + v_n^{(3)}.
\end{aligned}
$$

By (3.2),

$$
v_n^{(1)} = O(n^2 L_n^{-3}).
$$

The estimation of $v_n^{(2)}$ follows *mutatis mutandis* from that for $w_n^{(2)}$ by using (1.6) with $m = 6$:

$$
\begin{aligned}
v_n^{(2)} &\leq \sum_{1 \leq |\Delta| \leq K} \frac{\mathbb{E}\{(Y_{n,k} - \mu_{n,k})^6\}}{(\mu_{n,k_0} - \mu_{n,k_0+\Delta})^4} \\
&= O\left(n^2 L_n^{-7} \sum_{1 \leq |\Delta| \leq K} \frac{\Delta^6 e^{-3\Delta^2/L_n}}{(1 - e^{-\Delta^2/(3L_n)})^4}\right) \\
&= O\left(n^2 L_n^{-7/2} \int_{1/L_n}^{\infty} v^{-3/2}\,dv\right) \\
&= O(n^2 L_n^{-3}).
\end{aligned}
$$

For the last term $v_n^{(3)}$, we use again (2.3),

$$
\mathbb{V}\{Y_{n,k}\} \leq \mu_{n,k}^2 = O(L_n^{-1} r^{-2k} n^{2r})
$$

uniformly for $1 \leq k \leq n$, where $r > 0$ is any bounded real number. Substituting this into $v_n^{(3)}$ gives

$$
\begin{aligned}
v_n^{(3)} &= O\left(L_n^{-1} n^{2r} \sum_{|\Delta| \geq K} r^{-2k_0 - 2\Delta}\right) \\
&= O(L_n^{2/3} n^2 e^{-L_n^{1/3}})
\end{aligned}
$$

by taking $r = 1 + \text{sign}(\Delta) L_n^{-1/3}$.

This completes the proof of (3.4).

*Higher central moments of* $W_n$. The same analysis can be carried out for higher absolute central moments using (1.7). Then the same proof for $\mathbb{V}\{W_n\}$ gives (1.5) by using (1.7) with order $2s + 2$.

*Almost sure convergence.* We need first a tail bound for the width. By Markov's inequality (see page 160 of [29]; sometimes referred to as Chebyshev inequality),

$$\mathbb{P}\{|W_n - \mathbb{E}\{W_n\}| \geq \varepsilon\mathbb{E}\{W_n\}\} \leq \frac{\mathbb{E}\{|W_n - \mathbb{E}\{W_n\}|^s\}}{(\varepsilon\mathbb{E}\{W_n\})^s}$$
$$= O(\varepsilon^{-s} L_n^{-s})$$

for any $s > 0$ and $\varepsilon \in (0, 1)$.

From this estimate it follows, by applying the Borel–Cantelli lemma and by taking $s > 2$, that

$$\frac{W_{n_\ell}}{\mathbb{E}\{W_{n_\ell}\}} \longrightarrow 1 \qquad \text{almost surely,}$$

where $n_\ell := \lfloor e^{\sqrt{\ell}} \rfloor$, since $\sum_\ell L_{n_\ell}^{-s} = O(\sum_\ell \ell^{-s/2}) = O(1)$.

Now observe that

$$n_{\ell+1} - n_\ell = \Theta(n_\ell \ell^{-1/2}) = \Theta(n_\ell L_{n_\ell}^{-1}) = \Theta(\mathbb{E}\{W_{n_\ell}\} L_{n_\ell}^{-1/2}).$$

On the other hand, by construction, adding a new node to random recursive trees affects the value of $W_n$ by at most 1. Consequently,

$$\sup_{n_\ell \leq n < n_{\ell+1}} \max(|W_n - W_{n_\ell}|, |\mathbb{E}\{W_n\} - \mathbb{E}\{W_{n_\ell}\}|) \leq n_{\ell+1} - n_\ell$$
$$= \Theta(\mathbb{E}\{W_{n_\ell}\} L_{n_\ell}^{-1/2}).$$

So, deterministically,

$$\sup_{n_\ell \leq n < n_{\ell+1}} \left| \frac{W_n}{\mathbb{E}\{W_n\}} - \frac{W_{n_\ell}}{\mathbb{E}\{W_{n_\ell}\}} \right| = O\left( \frac{\mathbb{E}\{W_{n_\ell}\} L_{n_\ell}^{-1/2}}{\mathbb{E}\{W_{n_\ell}\} - (n_{\ell+1} - n_\ell)} \right)$$
$$= O(L_{n_\ell}^{-1/2})$$
$$= O(\ell^{-1/4}).$$

This completes the proof of (1.8).

*An alternative form to* (1.8). The same argument can be modified to show that

$$\frac{W_n}{n/\sqrt{2\pi L_n}} = 1 + O(L_n^{-1+\delta}) \tag{3.5}$$

almost surely for any fixed $\delta > 0$. The proof is modified from that for (1.8) as follows. By (1.3), we have

$$\frac{W_n}{n/\sqrt{2\pi L_n}} = \frac{W_n}{\mathbb{E}\{W_n\}}(1 + O(L_n^{-1})).$$

Instead of $n_\ell := \lfloor e^{\sqrt{\ell}} \rfloor$, we now take $n_\ell := \lfloor e^{\sqrt{\ell}/(2-\delta)} \rfloor$. Then, setting $\varepsilon = \varepsilon_n = L_n^{-1+\delta}$ in the proof, we deduce that, again by the Borel–Cantelli lemma,

$$\frac{W_{n_\ell}}{\mathbb{E}\{W_{n_\ell}\}} = 1 + O(\varepsilon_{n_\ell})$$

almost surely as $\ell \to \infty$ provided that $\varepsilon_{n_\ell}^{-s} L_{n_\ell}^{-s}$ is summable in $\ell$. This forces the choice $s > 2/\delta$. Next,

$$n_{\ell+1} - n_\ell = \Theta(\mathbb{E}\{W_{n_\ell}\}\ell^{-1/4}).$$

This proves (3.5).

*Almost sure convergence for $Y_{n,k}$.* We can also obtain strong convergence by the same argument for the profiles $Y_{n,k}$ in the central range $[L_n - L_n^{1-\varepsilon}, L_n + L_n^{1-\varepsilon}]$, where $\varepsilon \in (0,1)$. We prove that

$$(3.6) \qquad \sup_{L_n - L_n^{1-\varepsilon} \leq \kappa \leq L_n + L_n^{1-\varepsilon}} \left| \frac{Y_{n,\kappa}}{\mathbb{E}\{Y_{n,\kappa}\}} - 1 \right| \to 0$$

almost surely.

PROOF. Set $t_n := 2L_n^{1-\varepsilon}$ and $n_\ell := \lfloor e^{\sqrt{\ell}} \rfloor$. Using (1.7) and Markov's inequality used above, it is easy to see that

$$\sup_{L_{n_\ell} - t_{n_\ell} \leq \kappa \leq L_{n_\ell} + t_{n_\ell}} \left| \frac{Y_{n_\ell,\kappa}}{\mathbb{E}\{Y_{n_\ell,\kappa}\}} - 1 \right| \to 0,$$

almost surely as $\ell \to \infty$. By the union bound and the Borel–Cantelli lemma, this requires that we take $s$ so large that $L_{n_\ell}^{-s} t_{n_\ell}^{1+s}$ is summable in $\ell$. Any choice with $s > 3/\varepsilon - 1$ suffices for that purpose. Furthermore, by the monotonicity of $Y_{n,k}$ in $n$ for fixed $k$,

$$\sup_{n_\ell \leq n < n_{\ell+1}} \sup_{L_{n_\ell} - t_{n_\ell} \leq \kappa \leq L_{n_\ell} + t_{n_\ell}} |Y_{n,\kappa} - Y_{n_\ell,\kappa}|$$

$$\leq \sup_{L_{n_\ell} - t_{n_\ell} \leq \kappa \leq L_{n_\ell} + t_{n_\ell}} |Y_{n_{\ell+1},\kappa} - Y_{n_\ell,\kappa}|$$

and

$$\sup_{n_\ell \leq n < n_{\ell+1}} \sup_{L_{n_\ell} - t_{n_\ell} \leq \kappa \leq L_{n_\ell} + t_{n_\ell}} |\mathbb{E}\{Y_{n,\kappa}\} - \mathbb{E}\{Y_{n_\ell,\kappa}\}|$$

$$\leq \sup_{L_{n_\ell} - t_{n_\ell} \leq \kappa \leq L_{n_\ell} + t_{n_\ell}} |\mathbb{E}\{Y_{n_{\ell+1},\kappa}\} - \mathbb{E}\{Y_{n_\ell,\kappa}\}|.$$

Thus,

$$\sup_{n_\ell \leq n < n_{\ell+1}} \left| \frac{Y_{n,\kappa}}{\mathbb{E}\{Y_{n,\kappa}\}} - \frac{Y_{n_\ell,\kappa}}{\mathbb{E}\{Y_{n_\ell,\kappa}\}} \right|$$

$$\leq \frac{Y_{n_{\ell+1},\kappa}}{\mathbb{E}\{Y_{n_\ell,\kappa}\}} - \frac{Y_{n_\ell,\kappa}}{\mathbb{E}\{Y_{n_\ell,\kappa}\}}$$

$$\leq \frac{Y_{n_{\ell+1},\kappa}}{\mathbb{E}\{Y_{n_{\ell+1},\kappa}\}} - \frac{Y_{n_\ell,\kappa}}{\mathbb{E}\{Y_{n_\ell,\kappa}\}} + \left( \frac{\mathbb{E}\{Y_{n_{\ell+1},\kappa}\}}{\mathbb{E}\{Y_{n_\ell,\kappa}\}} - 1 \right) \frac{Y_{n_{\ell+1},\kappa}}{\mathbb{E}\{Y_{n_{\ell+1},\kappa}\}}.$$

Putting the supremum over $L_{n_\ell} - t_{n_\ell} \leq \kappa \leq L_{n_\ell} + t_{n_\ell}$ in front of all of the latter inequalities, we see that both terms tend to zero almost surely provided that

$$\lim_{\ell \to \infty} \sup_{L_{n_\ell} - t_{n_\ell} \leq \kappa \leq L_{n_\ell} + t_{n_\ell}} \left| \frac{\mathbb{E}\{Y_{n_{\ell+1},\kappa}\}}{\mathbb{E}\{Y_{n_\ell,\kappa}\}} - 1 \right| = 0.$$

This follows from an extension of the Taylor series estimate used in (1.1); indeed, the estimate (see [25])

$$\mu_{n,k} = \frac{L_n^k}{\Gamma(1 + k/L_n)k!}(1 + O(L_n^{-1})) \qquad (k = O(L_n))$$

is sufficient for our use.

Thus we have shown that

$$\sup_{n_\ell \leq n < n_{\ell+1}} \sup_{L_{n_\ell} - t_{n_\ell} \leq \kappa \leq L_{n_\ell} + t_{n_\ell}} \left| \frac{Y_{n,\kappa}}{\mathbb{E}\{Y_{n,\kappa}\}} - 1 \right| \to 0$$

almost surely. An additional argument shows that for $\ell$ large enough, $[L_n - L_n^{1-\varepsilon}, L_n + L_n^{1-\varepsilon}]$ is contained in $[L_{n_\ell} - t_{n_\ell}, L_{n_\ell} + t_{n_\ell}]$ for $n_\ell \leq n < n_{\ell+1}$, thus concluding the proof of (3.6).  □

3.3. *Level reaching the width.* We now prove Theorem 1.3. For $|\Delta| > |\hat{k} - k_0|$ and $B > 1$,

$$\mathbb{P}(k^* = k_0 + \Delta) = \mathbb{P}(W_n = Y_{n,k_0+\Delta})$$

$$\leq \mathbb{P}(Y_{n,k_0+\Delta} > Y_{n,k_0})$$

$$= \mathbb{P}(Y_{n,k_0+\Delta} - \mu_{n,k_0+\Delta} > Y_{n,k_0} - \mu_{n,k_0} + \mu_{n,k_0} - \mu_{n,k_0+\Delta})$$

$$\leq \mathbb{P}\left( Y_{n,k_0+\Delta} - \mu_{n,k_0+\Delta} \geq \frac{1}{2}(\mu_{n,k_0} - \mu_{n,k_0+\Delta}) \right)$$

$$+ \mathbb{P}\left( Y_{n,k_0} - \mu_{n,k_0} \leq -\frac{1}{2}(\mu_{n,k_0} - \mu_{n,k_0+\Delta}) \right)$$

$$\leq \frac{2^B \mathbb{E}|Y_{n,k_0+\Delta} - \mu_{n,k_0+\Delta}|^B}{(\mu_{n,k_0} - \mu_{n,k_0+\Delta})^B} + \frac{2^B \mathbb{E}|Y_{n,k_0} - \mu_{n,k_0}|^B}{(\mu_{n,k_0} - \mu_{n,k_0+\Delta})^B}$$

by Markov's inequality. By (1.1), we obtain an estimate similar to (3.1) for $\mu_{n,k_0} - \mu_{n,k_0+\Delta}$, which together with (3.2) gives

$$\mathbb{P}(k^* = k_0 + \Delta) = O\left(\frac{\Delta^B L_n^{-B} e^{-B\Delta^2/(2L_n)}}{(1 - e^{-\Delta^2/(3L_n)})^B} + \frac{L_n^{-B}}{(1 - e^{-\Delta^2/(3L_n)})^B}\right)$$

$$= O(\Delta^{-B} + \Delta^{-2B} \vee L_n^{-B})$$

$$= O(\Delta^{-B})$$

uniformly for $1 \leq |\Delta| \leq K$. It follows that there exists a $T_0 > 1$ such that for $T > T_0$,

$$\mathbb{P}(|k^* - k_0| \geq T) = O\left(\sum_{T \leq |\Delta| \leq K} \Delta^{-B}\right) + \mathbb{P}(|k^* - k_0| \geq K)$$

$$= O(T^{1-B}) + \mathbb{P}(|k^* - k_0| \geq K).$$

The tail probability $\mathbb{P}(|k^* - k_0| \geq K)$ is estimated as follows, where we let $k_1 := \lfloor \sqrt{L_n} \rfloor$:

$$\mathbb{P}(|k^* - k_0| \geq K) \leq \mathbb{P}\left(\max_{|k-k_0| \geq K} Y_{n,k} \geq Y_{n,k_0}\right)$$

$$\leq \mathbb{P}\left(\max_{|k-k_0| \geq K} Y_{n,k} \geq \mu_{n,k_0+k_1}\right) + \mathbb{P}(Y_{n,k_0} < \mu_{n,k_0+k_1})$$

$$\leq \mu_{n,k_0+k_1}^{-1} \sum_{|k-k_0| \geq K} \mu_{n,k} + \frac{\mathbb{V}\{Y_{n,k_0}\}}{(\mu_{n,k_0} - \mu_{n,k_0+k_1})^2}$$

$$= O(L_n^{1/3} e^{-L_n^{1/3}/2} + L_n^{-2}),$$

which tends to zero as $n \to \infty$, where we used again (2.3) to bound $\sum_{|k-k_0| \geq K} \mu_{n,k}$. Since $B > 1$ is arbitrary, this proves Theorem 1.3.

*Limit distribution of $W_n$?* It is known that the centered and normalized random variables $(Y_{n,k} - \mu_{n,k})/\sqrt{\mathbb{V}\{Y_{n,k}\}}$ do not converge to a fixed limit law when $k = L_n + O(1)$ due to periodicity; see [21]. The origin of the periodicity lies at the second-order term in the asymptotic expansion of $\mu_{n,L_n+O(1)}$,

$$\mu_{n,k_0+\ell} = \frac{n}{\sqrt{2\pi L_n}}\left(1 + \frac{p_\ell(\{L_n\})}{L_n} + O(L_n^{-2})\right) \qquad (\ell \in \mathbb{Z}),$$

where $\{x\}$ denotes the fractional part of $x$ and

$$p_\ell(x) := -\frac{1}{2}\left(x - \ell - \frac{3}{2} + \gamma\right)^2 - \frac{\gamma}{2} + \frac{\pi^2}{12} + \frac{1}{24}.$$

This periodic second-order term is the origin of all fluctuations of higher central moments. Note that

$$(3.7) \qquad \max_{\substack{\ell \in \mathbb{Z} \\ x \in [0,1]}} p_\ell(x) = \begin{cases} p_{-1}(x), & \text{if } x \in [0, 1-\gamma], \\ p_0(x), & \text{if } x \in [1-\gamma, 1]; \end{cases}$$

compare (1.4).

The main open question is the limit distribution (if it exists) of $W_n$. Simulations seem to indicate the closeness of the histogram of $W_n$ to that of $Y_{n,k_0-1}$ when $\{L_n\} < 1 - \gamma$ and to that of $Y_{n,k_0}$ when $\{L_n\} > 1 - \gamma$; see Figure 1.

**4. Width of general random log trees.** Our methods of proof for recursive trees can be formulated in a general, simple framework described below, which gives sufficient conditions we need to obtain asymptotics of the width of general random log trees. We first describe the estimates we need to handle the width of general random log trees and then discuss a few concrete examples in the remaining sections.

4.1. *An analytic scheme for the profile and width.* We start from a general framework for the moments of the profiles of random log trees. The same notations as those for recursive trees are used in this section. So we denote by $Y_{n,k}$ the profile of the random log trees in question, the initial condition being immaterial for our purpose. We impose the following three conditions.

CONDITION I (*Quasi-power form for the expected profile polynomial*). The generating polynomial of the expected profile $\mu_{n,k} := \mathbb{E}\{Y_{n,k}\}$ of the random log trees in question satisfies asymptotically

$$(4.1) \qquad \Xi_n(u) := \sum_k \mu_{n,k} u^k = g(u) n^{f(u)} (1 + O(n^{-\varepsilon}))$$

uniformly for $|u - 1| \leq \varepsilon_0$, $u \in \mathbb{C}$, $\varepsilon_0 > 0$. Here $g$ and $f$ are analytic functions in $|u - 1| \leq \varepsilon_0$ and satisfy $g(1) = f(1) = 1$.

CONDITION II (*Regularity condition for the expected profile polynomial*). The estimate

$$(4.2) \qquad |\Xi_n(u)| = O(n^{1-\varepsilon})$$

holds uniformly for $\{u \in \mathbb{C} : 1 - \varepsilon_1 \leq |u| \leq 1 + \varepsilon_1\} \setminus \{u \in \mathbb{C} : |u - 1| \leq \varepsilon_0\}$, where $0 < \varepsilon_1 < \varepsilon_0$; see Figure 2 for a plot.

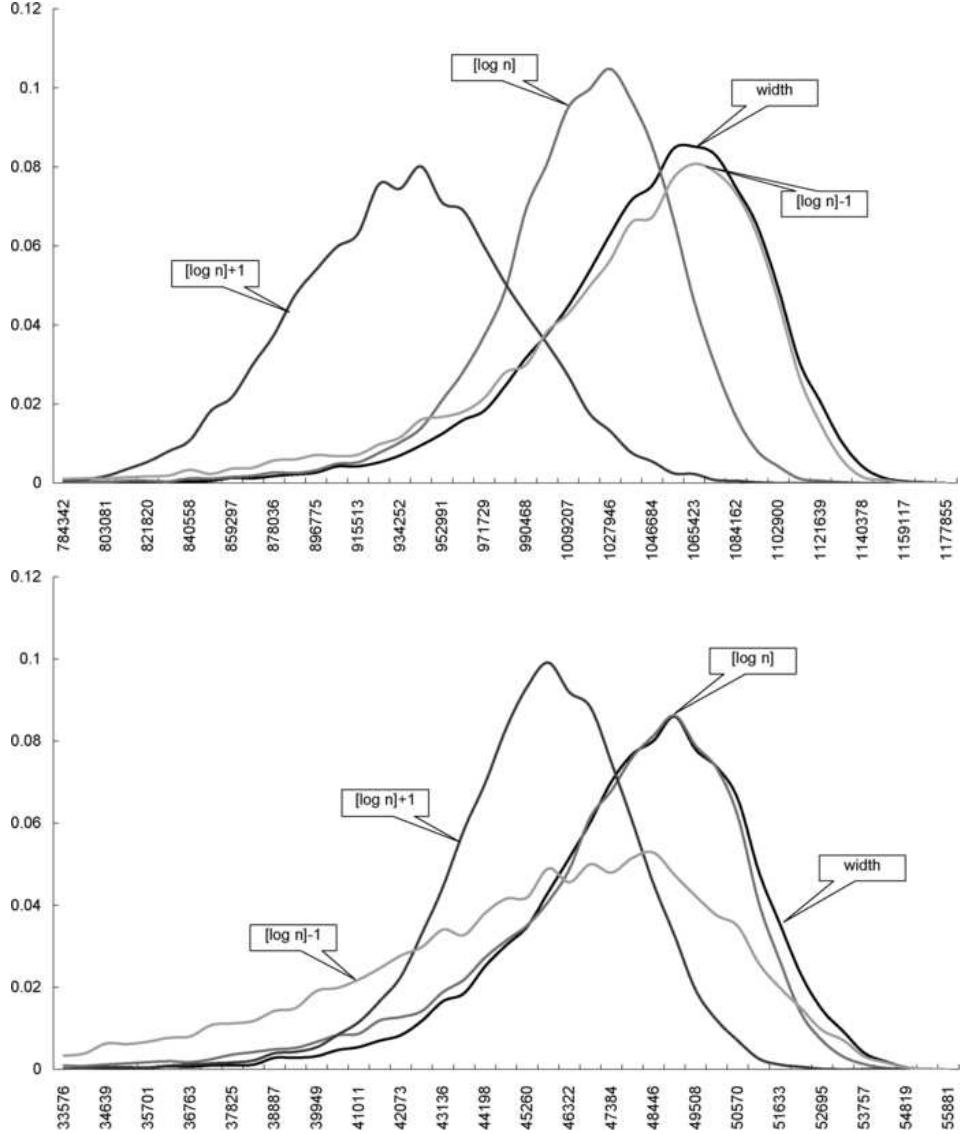

Fig. 1. *Simulated histograms for $W_n$ and $Y_{n,\lfloor L_n \rfloor + \ell}$ for $\ell = -1, 0, 1$, where $n = 10^7$ (top), for which the expected width is asymptotically reached at $\lfloor L_n \rfloor - 1$ ($\{L_n\} < 1 - \gamma$), and $n = 404960$ (bottom), for which the expected width is asymptotically reached at $\lfloor L_n \rfloor$ ($\{L_n\} > 1 - \gamma$); see (1.4) and (3.7).*

CONDITION III (*Asymptotic estimates for the central moments*). The central moments of $Y_{n,k}$ satisfy

$$(4.3) \qquad \mathbb{E}\{(Y_{n,k} - \mu_{n,k})^m\} = O(|\Delta|^m L_n^{-m} \mu_{n,k}^m) \qquad (m \geq 0)$$

uniformly for $|\Delta| = o(L_n)$.

THEOREM 4.1 (Width regularity of general random log trees). *Under Conditions* I, II *and* III, *the width of the random log trees in question satisfies the properties*

$$\mathbb{E}\{W_n\} = \frac{n}{\sqrt{2\pi\sigma^2 L_n}}(1 + O(L_n^{-1})),$$

(4.4)        $$\mathbb{E}\{|W_n - \mathbb{E}\{W_n\}|^s\} = O(n^s L_n^{-3s/2}) \qquad (s \geq 0),$$

$$\mathbb{P}(|k^* - f'(1)L_n| \geq T) = O(T^{-B}).$$

*The last estimate holds for every $B > 0$ and $T > T_0$ for some $T_0 > 1$. Furthermore, if inserting a new node to the tree changes the width by at most a bounded quantity, then we also have*

(4.5)        $$\frac{W_n}{\mathbb{E}\{W_n\}} \longrightarrow 1 \qquad almost\ surely.$$

For ease of reference, we will refer to the properties (4.4) by saying that the random log trees are width-regular with parameters $(f'(1), \sigma^2)$.

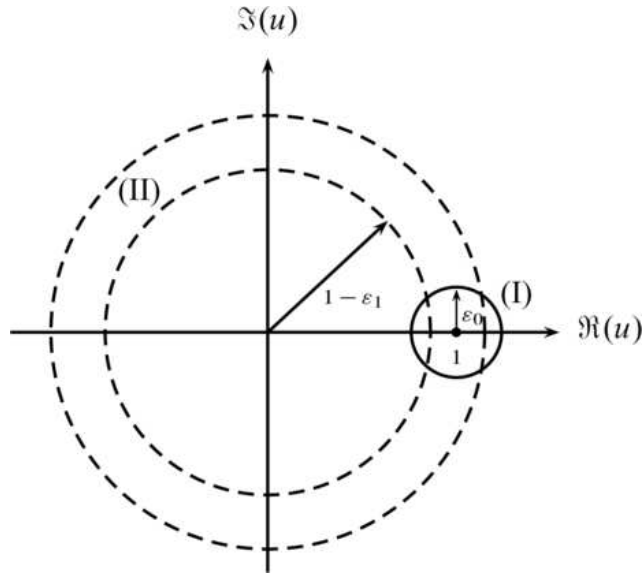

FIG. 2. *The different regions in the complex u-plane as used by Conditions* I *and* II *for Theorem* 4.1.

PROOF (Sketch). First, Conditions I and II imply, by standard application of the saddle point method, that

$$(4.6) \qquad \mu_{n,k} = \frac{n}{\sqrt{2\pi\sigma^2 L_n}} e^{-\Delta^2/(2\sigma^2 L_n) + O(|\Delta|^3/L_n^2)} \left(1 + O\left(\frac{1+|\Delta|}{L_n}\right)\right)$$

uniformly for $|\Delta| \le L_n^{2/3}$, where $\Delta := k - f'(1)L_n$ and

$$\sigma = \sqrt{f'(1) + f''(1)}.$$

Note that to prove the estimate (4.6), we used (4.1) and (4.2) only when $u = e^{i\theta}$, $\theta \in \mathbb{R}$. However, the uniform estimates (4.1) and (4.2) in a complex neighborhood of unity also yield, by Cauchy's integral representation,

$$(4.7) \qquad \mu_{n,k} = O(L_n^{-1/2} r^{-k} n^{f(r)} + r^{-k} n^{1-\varepsilon}) = O(L_n^{-1/2} r^{-k} n^{f(r)})$$

uniformly for all $k = 0, \ldots, n$, where $r = 1 + o(1)$. This crude estimate is sufficient for our purpose in bounding all error terms involved.

The remaining proofs follow closely those used for recursive trees, details being omitted here. □

Note that our proof for the almost sure convergence (4.5) requires the estimate for $\mathbb{E}\{|W_n - \mathbb{E}\{W_n\}|^{2+\varepsilon}\}$ for which (4.3) with $m = 8$ suffices. Similarly, the estimate for the expected width needs (4.3) with $m = 4$.

4.2. *The moments estimates (4.3).* Theorem 4.1 reduces the proofs of the width-regularity properties to those for the three hard estimates (4.1), (4.2) and (4.3). While the two estimates (4.1) and (4.2) are more tree dependent and may rely on different analytic tools, we indicate in this subsection how the moments estimates (4.3) can in many cases be deduced from (4.1) and (4.2), coupling with suitable "asymptotic transfer" for the underlying bivariate recurrence.

*Recurrence of profile.* The profiles of many random log trees that arise from data structures, analysis of algorithms and discrete probability are of the form (see [2, 8, 28, 32])

$$(4.8) \qquad Y_{n,k} \stackrel{d}{=} \sum_{1 \le j \le h} Y_{I_{n,j},k-1}^{(j)} \qquad (n \ge 2; k \ge 1),$$

with $Y_{n,0} = 1$ for $n \ge 1$, where $h \ge 2$, the $Y_{n,k}^{(j)}$'s are independent copies of $Y_{n,k}$ and the underlying splitting distribution satisfies $\sum_{1 \le j \le h} I_{n,j} = n - \kappa$ for some integer $\kappa \ge 0$. Physically, the root of such random log trees has at most $h$ subtrees, each of which has the same (recursive) structure; the distribution of the size of the $j$th subtree of the root is described by $I_{n,j}^{(j)}$ and $\kappa$ represents the number of nodes retained at the root.

Then the moments of $Y_{n,k}$ satisfy a recurrence of the form

$$(4.9) \qquad a_{n,k} = h \sum_{0 \le j < n} \pi_{n,j} a_{j,k-1} + b_{n,k},$$

where $\pi_{n,j} = \mathbb{P}(I_{n,1} = j)$ satisfies $\sum_j \pi_{n,j} = 1$ and the $b_{n,k}$'s are known. For our purpose, we can always assume that $b_{n,k} = 0$ for $k < 0$ and $k \ge n$.

*Higher central moments.* If $Y_{n,k}$ satisfies the distributional recurrence (4.8), then the central moments $P_{n,k}^{(m)} := \mathbb{E}\{(Y_{n,k} - \mu_{n,k})^m\}$ can be recursively computed by the recurrence

$$P_{n,k}^{(m)} = h \sum_{0 \le j < n} \pi_{n,j} P_{j,k-1}^{(m)} + Q_{n,k}^{(m)},$$

where

$$Q_{n,k}^{(m)} := \sum_{\substack{i_0 + i_1 + \cdots + i_h = m \\ 0 \le i_1, \dots, i_h < m \\ 0 \le i_0 \le m}} \sum_{j_1 + \cdots + j_h = n - \kappa} \mathbb{P}\{I_{n,1} = j_1, \dots, I_{n,h} = j_h\}$$

$$\times P_{j_1,k-1}^{(i_1)} \cdots P_{j_h,k-1}^{(i_h)} \nabla_{n,k}^{i_0}(\mathbf{j}) \qquad (m \ge 2),$$

with $[\,\mathbf{j} = (j_1, \dots, j_h)\,]$

$$(4.10) \qquad \nabla_{n,k}(\mathbf{j}) := \sum_{1 \le \ell \le h} \mu_{j_\ell,k-1} - \mu_{n,k}.$$

By Cauchy's integral formula and (4.1), we have

$$\nabla_{n,k}(\mathbf{j}) = \frac{1}{2\pi i} \int_{\substack{|u|=r \\ |u-1| \le \varepsilon}} g(u) u^{-k-1} n^{f(u)} \varphi\left(u; \frac{\mathbf{j}}{n}\right)$$

$$\times \left(1 + O\left(\sum_{1 \le \ell \le h} \frac{1}{(j_\ell + 1)^\varepsilon}\right)\right) du + O(r^{-k} n^{1-\varepsilon}),$$

where

$$\varphi\left(u; \frac{\mathbf{j}}{n}\right) := \sum_{1 \le \ell \le h} u\left(\frac{j_\ell}{n}\right)^{f(u)} - 1.$$

Since $\sum_{1 \le \ell \le h} j_\ell = n + O(1)$, we deduce, by expanding $\varphi(u; x)$ at $u = 1$, the two estimates

$$(4.11) \quad \nabla_{n,k}(\mathbf{j}) = \begin{cases} O((|r-1| + L_n^{-1/2}) L_n^{-1/2} r^{-k} n^{f(r)}) & (r = 1 + o(1)), \\ O(|\Delta| L_n^{-3/2} n), \end{cases}$$

where the first estimate holds uniformly for all tuples $(j_1, \ldots, j_h)$ and $1 \leq k < n$, and the second holds for all tuples $(j_1, \ldots, j_h)$ and $\Delta = O(\sqrt{L_n})$. Note that if we take $r$ to be the solution near unity of the equation $rf'(r) = k/L_n = f'(1) + \Delta/L_n$, then $r = 1 + \Delta/(\sigma^2 L_n) + O(\Delta^2/L_n^2)$ and

$$(4.12) \qquad r^{-k} n^{f(r)} = n e^{-\Delta^2/(2\sigma^2 L_n) + O(|\Delta|^3 L_n^{-2})}$$

uniformly for $\Delta = O(L_n^{2/3})$. This means that the first estimate in (4.11) is not tight when $\Delta = o(\sqrt{L_n})$, which is the reason why we need the second estimate.

Once the precise estimates for $\nabla_{n,k}(\mathbf{j})$ are available, the remaining proof for (4.3) is then reduced to the derivation of suitable "asymptotic transfer" for the recurrence (4.9) via which one deduces estimate for $a_{n,k}$ from that for $b_{n,k}$. Instead of further abstraction for general random log trees, we will give more details for specific trees below.

## 5. Random quad trees and grid trees.
We start from quad trees, which are useful data structures for spatial points, and then indicate the estimates needed for the more general grid trees proposed in [8]. We show that both classes of trees are width-regular.

*Random quad trees and their construction.* Quad trees were proposed by Finkel and Bentley [15]. The first probabilistic analysis of the typical depth of a node, the expected profile, the height, and the partial match cost was carried out by Flajolet, Gonnet, Puech and Robson [16, 17] (work carried out in 1988), Devroye and Laforest [9] and Flajolet, Labelle, Laforest and Salvy [18].

Given a sequence of $n$ points independently and uniformly chosen from $[0,1]^d$, the random (point) quad tree associated with this random sample is constructed by placing the first point at the root, which splits the space into $2^d$ hyperrectangles, each corresponding to one of the $2^d$ subtrees of the root. Points that fall in each hyperrectangle are directed to the corresponding subtree and are constructed recursively. For more information on quad trees, see [18, 28, 32, 41] and the references therein.

*The profile.* By such a construction, the profile $Y_{n,k}$ satisfies (4.8) with $h = 2^d$, $\kappa = 1$ and

$$\pi_{n,\mathbf{j}} := \mathbb{P}(I_{n,1} = j_1, \ldots, I_{n,2^d} = j_{2^d})$$

$$= \binom{n-1}{j_1, \ldots, j_{2^d}}$$

$$\times \int_{[0,1]^d} \prod_{\substack{1 \leq \ell \leq 2^d \\ \ell-1=(b_1,\ldots,b_d)_2}} \left( \prod_{1 \leq i \leq d} b_i(1-x_i) + (1-b_i)x_i \right)^{j_\ell} d\mathbf{x},$$

where $(b_1, \ldots, b_d)_2$ denotes the binary representation of $\ell - 1$ [prefixed by zeros if $\lfloor \log_2(\ell - 1) \rfloor < d - 1$] and $d\mathbf{x} = dx_1 \cdots dx_d$.

*The underlying recurrence.* From the expression for $\pi_{n,\mathbf{j}}$, it follows that all moments of $Y_{n,k}$ satisfy (4.9) with (see [18])

$$(5.1) \qquad \pi_{n,j} = \frac{1}{n} \sum_{j < j_1 \leq \cdots \leq j_{d-1} \leq n} \frac{1}{j_1 \cdots j_{d-1}} \qquad (0 \leq j < n).$$

In particular, the expected profile $\mu_{n,k}$ satisfies the estimates (4.1) and (4.2) with $f(u) = 2u^{1/d} - 1$ and

$$g(u) := \frac{1}{\Gamma(2u^{1/d})^d (2u^{1/d} - 1)}$$

$$\times \prod_{1 \leq \ell < d} \frac{\Gamma(2u^{1/d}(1 - e^{2\ell\pi i/d}))}{\Gamma(2 - 2u^{1/d} e^{2\ell\pi i/d})};$$

see [4] and [20]. The exact form of $g$ is less important for our purpose; the analyticity of $g$ for $u$ near unity is, however, technically useful. Note that

$$f'(1) = \frac{2}{d}, \qquad \sigma^2 = \frac{2}{d^2}.$$

*Recurrence of $P_{n,k}^{(m)} := \mathbb{E}\{(Y_{n,k} - \mu_{n,k})^m\}$.* Obviously, $P_{n,k}^{(0)} = 1$, $P_{n,k}^{(1)} = 0$ and $P_{n,k}^{(m)}$ satisfies the recurrence

$$P_{n,k}^{(m)} = 2^d \sum_{0 \leq j < n} \pi_{n,j} P_{j,k-1}^{(m)} + Q_{n,k}^{(m)} \qquad (m \geq 2),$$

where

$$Q_{n,k}^{(m)} := \sum_{(i_0, \ldots, i_{2d}) \in \mathcal{I}_m} \binom{m}{i_0, \ldots, i_{2d}}$$

$$\times \sum_{j_1 + \cdots + j_{2d} = n-1} \pi_{n,\mathbf{j}} P_{j_1, k-1}^{(i_1)} \cdots P_{j_{2d}, k-1}^{(i_{2d})} \nabla_{n,k}^{i_0}(\mathbf{j}).$$

Here $\nabla_{n,k}(\mathbf{j})$ is given in (4.10) with $h = 2^d$ there and

$$\mathcal{I}_m := \{(i_0, \ldots, i_{2d}) \in \{0, \ldots, m\} \times \{0, \ldots, m-1\}^d : i_0 + \cdots + i_{2d} = m\}.$$

Following the proof pattern for recursive trees and the discussions in Section 4.2, we prove, based on the estimates (4.11), the two bounds

$$(5.2) \quad P_{n,k}^{(m)} = \begin{cases} O((|r-1|^m + L_n^{-m/2})L_n^{-m/2} r^{-mk} n^{mf(r)}), & r = 1 + o(1), \\ O(|\Delta|^m L_n^{-3m/2} n^m), \end{cases}$$

the first being uniform for $1 \leq k < n$ and the second uniform for $\Delta := k - f'(1)L_n = O(\sqrt{L_n})$. These two estimates imply (4.3) by (4.7) and (4.12).

*An asymptotic transfer for the double-indexed recurrence.* To justify the width-regularity properties (4.4), it remains to prove the two estimates in (5.2). For an exact solution for (4.9) similar to (2.9), see [18]. Here we use a different inductive argument, which is easily amended for other varieties of trees.

LEMMA 1. *Assume that $a_{n,k}$ satisfies* (4.9) *with $\pi_{n,j}$ given in* (5.1) *and that $a_{n,0}, a_{1,k} = O(1)$. If*

$$|b_{n,k}| \leq c|k - f'(1)L_n|^\lambda L_n^\beta \rho^{-k} n^\alpha$$

*for $n \geq 1$ and $1 \leq k \leq n$, where $\lambda \geq 0$, $\beta \in \mathbb{R}$, $c > 0$ and the two real numbers $\alpha, \rho > 0$ satisfy $\rho < ((\alpha+1)/2)^d$, then*

$$(5.3) \qquad |a_{n,k}| \leq C_0|k - f'(1)L_n|^\lambda L_n^\beta \rho^{-k} n^\alpha$$

*for $n \geq 1$ and $1 \leq k \leq n$, where $C_0 > 0$ is chosen so large that $C_0 \geq c/(1 - \varepsilon - 2^d \rho/(\alpha+1)^d)$.*

PROOF. We apply induction on $k$ and $n$. The boundary conditions are easily checked by taking $C_0$ sufficiently large. We may assume that $|k - f'(1)L_n| \to \infty$, for otherwise we need only modify the value of $c$. By the induction hypothesis, we have (see [4])

$$|a_{n,k}| \leq c|k - f'(1)L_n|^\lambda L_n^\beta \rho^{-k} n^\alpha$$

$$+ 2^d C_0 \rho^{1-k} n^{-1} \sum_{1 \leq j < j_1 \leq \cdots \leq j_{d-1} \leq n} \frac{|k - 1 - f'(1)L_j|^\lambda L_j^\beta j^\alpha}{j_1 \cdots j_{d-1}}$$

$$\sim c|k - f'(1)L_n|^\lambda L_n^\beta \rho^{-k} n^\alpha$$

$$+ \frac{2^d C_0}{(d-1)!} \rho^{1-k} n^{-1} \sum_{1 \leq j < n} |k - f'(1)L_j|^\lambda L_j^\beta j^\alpha \left(\log \frac{n}{j}\right)^{d-1}$$

$$= cn^\alpha L_n^\beta \rho^{-k} + \frac{2^d}{(\alpha+1)^d} C_0 (1 + o(1))|k - f'(1)L_n|^\lambda L_n^\beta \rho^{1-k} n^\alpha;$$

thus (5.3) follows by properly tuning $C_0$ [since $\rho < ((\alpha+1)/2)^d$]. □

*Asymptotics of $P_{n,k}^{(m)}$.* We prove by induction the first bound in (5.2). Assume $m \geq 2$.

Consider $Q_{n,k}^{(m)}$. As for recursive trees, we distinguish between two cases. If $j_1, \ldots, j_{2^d} \geq n/L_n^m$, then $L_{j_\ell} \sim L_n$ for $\ell = 1, \ldots, 2^d$ and we have

$$\sum_{(i_0, \ldots, i_{2^d}) \in \mathcal{I}_m} \binom{m}{i_0, \ldots, i_{2^d}}$$

$$(5.4) \qquad\qquad \times \sum_{n/L_n^m \leq j_1, \ldots, j_{2^d} < n} \pi_{n,\mathbf{j}} P_{j_1,k-1}^{(i_1)} \cdots P_{j_{2^d},k-1}^{(i_{2^d})} \nabla_{n,k}^{i_0}(\mathbf{j})$$

$$= O((|r-1|^m + L_n^{-m/2}) L_n^{-m/2} r^{-mk} n^{mf(r)}).$$

We now assume that one of the $j_\ell$'s, say $j_1$, is less than $n/L_n^m$. We may furthermore assume that the corresponding index $i_1$ of $j_1$ is nonzero; for otherwise, if all $i_\ell = 0$ for those $j_\ell$'s with $j_\ell \leq n/L_n^m$, then the bound on the right-hand side of (5.4) obviously holds, since all other $j_\ell$'s satisfy $L_{j_\ell} \sim L_n$. Terms in $Q_{n,k}^{(m)}$ with $i_1 \geq 1$ and $j_1 \leq n/L_n^m$ are bounded above by

$$O\left(r^{-mk} \sum_{(i_0, \ldots, i_{2^d}) \in \mathcal{I}_m} n^{(m-i_1)f(r)} \sum_{j_1 \leq n/L_n^m} \pi_{n,j_1} j_1^{i_1 f(r)}\right) = O(L_n^{-m} r^{-mk} n^{mf(r)}).$$

This proves that

$$Q_{n,k}^{(m)} = O((|r-1|^m + L_n^{-m/2}) L_n^{-m/2} r^{-mk} n^{mf(r)}).$$

Thus the first estimate in (5.2) holds by applying the $O$-transfer of Lemma 1. The proof of the second estimate in (5.2) follows by applying the same inductive argument; the details are omitted here.

Consequently, *random quad trees are width-regular with parameters* $(2/d, 2/d^2)$ *and with the almost sure convergence* (4.5). All our results are new except when $d = 1$, for which quad trees reduce to binary search trees and the almost sure convergence in (4.4) was derived in [3], and the expected width in [13] (with a weaker error term).

*Random grid trees.* Grid trees were first proposed by Devroye [8] and represent one of the extensions of quad trees. Instead of placing the first element in the given sequence at the root (as in quad trees), we fix an integer $m \geq 2$ and place the first $m-1$ elements at the root, which then split the space into $m^d$ hyperrectangles (called grids). The remaining construction is similar to that for quad trees.

In this case, we have $h = m^d$ and $(j_0 := j, j_d := n - m + 1)$

$$\pi_{n,j} = \sum_{j \leq j_1 \leq \cdots \leq j_{d-1} \leq n-m+1} \prod_{1 \leq \ell \leq d} \frac{\binom{j_\ell - j_{\ell-1} + m - 2}{m-2}}{\binom{j_\ell + m - 1}{m-1}},$$

and (4.1) and (4.2) hold by applying the approach proposed in [4], where $f(u)$ satisfies

$$((f(u)+1) \cdots (f(u)+m-1))^d = m!^d u \qquad (m \geq 2; d \geq 1),$$

with $f(1) = 1$. An $O$-transfer similar to that given in Lemma 1 can also be derived by noting that

$$\sum_{1 \leq j < n} \pi_{n,j} |k - 1 - f'(1)L_j|^\lambda L_j^\beta j^\alpha$$

$$= \sum_{1 \leq j \leq j_1 \leq \cdots \leq j_{d-1} \leq n-m+1} |k - 1 - f'(1)L_j|^\lambda L_j^\beta j^\alpha \prod_{1 \leq \ell \leq d} \frac{\binom{j_\ell - j_{\ell-1} + m - 2}{m-2}}{\binom{j_\ell + m - 1}{m-1}}$$

$$\sim \frac{(m-1)^d}{n}$$

$$\times \sum_{1 \leq j_{d-1} \leq n} \frac{1}{j_{d-1}} \left(1 - \frac{j_{d-1}}{n}\right)^{m-2}$$

$$\times \sum_{1 \leq j_{d-2} \leq j_{d-1}} \cdots \sum_{1 \leq j_1 \leq j_2} \frac{1}{j_1} \left(1 - \frac{j_1}{j_2}\right)^{m-2}$$

$$\times \sum_{1 \leq j \leq j_1} |k - f'(1)L_j|^\lambda L_j^\beta j^\alpha \left(1 - \frac{j}{j_1}\right)^{m-2}$$

$$\sim (m-1)^d \left(\frac{\Gamma(m-1)\Gamma(\alpha+1)}{\Gamma(m+\alpha)}\right)^d |k - f'(1)L_n|^\lambda L_n^\beta n^\alpha,$$

so that the same type of asymptotic transfer there holds when $\alpha, \rho > 0$ satisfy the inequality

$$\rho < \left(\frac{(\alpha+1)\cdots(\alpha+m-1)}{m!}\right)^d.$$

Thus *random grid trees are width-regular with*

$$f'(1) = \frac{1}{d(H_m - 1)}, \qquad \sigma^2 = \frac{H_m^{(2)} - 1}{d^2(H_m - 1)^3}$$

*and with the almost sure convergence* (4.5), *where* $H_m := \sum_{1 \leq j \leq m} 1/j$ *and* $H_m^{(2)} := \sum_{1 \leq j \leq m} 1/j^2$. Note that $d = 1$ corresponds to $m$-ary search trees (see [32]), and $m = 2$ corresponds to quad trees. No martingale structure is known for grid trees for general $(m, d)$. Our results are new.

## 6. Generalized $m$-ary search trees.

The $m$-ary search tree, proposed by Muntz and Uzgalis [38], generalizes the binary search tree. For the random version built from a random permutation of $\{1, \ldots, n\}$, early results on the typical depth and expected profile are from [33, 34]. These trees in turn led to the generalized $m$-ary search trees of Hennequin [23] (see also [5]). Instead of placing the first $m - 1$ elements in the given sequence of numbers at the root (as in $m$-ary search trees), we choose a random sample of $m(t+1)-1$ elements, where $m \geq 2$ and $t \geq 0$, and sort it in increasing order. Then use the $(t+1)$st, the $2(t+1)$st, $\ldots$ and the $(m-1)(t+1)$st smallest elements in the sample to partition the original sample into $m$ groups, corresponding to the $m$ subtrees of the root. Elements that fall in each subtree are constructed recursively in the same way and the process stops as long as the subtree size is less than $m(t+1)-1$, which can then be arranged arbitrarily, since asymptotically this will have a limited effect.

In this case, the profile $Y_{n,k}$ satisfies (4.8) with $h = m$ and

$$\mathbb{P}(I_{n,1} = j_1, \ldots, I_{n,m} = j_m) = \frac{\binom{j_1}{t} \cdots \binom{j_m}{t}}{\binom{n}{m(t+1)-1}}.$$

Furthermore, (4.1) and (4.2) hold with $f(u)$ satisfying the equation (see [32])

$$(f(u) + t + 1) \cdots (f(u) + m(t+1) - 1) = \frac{(m(t+1))!}{(t+1)!} u,$$

with $f(1) = 1$, where $m \geq 2$ and $t \geq 0$; see [5, 6, 7] for the asymptotic tools needed (based on differential equations). Straightforward computation gives

$$f'(1) = \frac{1}{H_{m(t+1)} - H_{t+1}}, \qquad \sigma^2 = \frac{H_{m(t+1)}^{(2)} - H_{t+1}^{(2)}}{(H_{m(t+1)} - H_{t+1})^3}.$$

The estimate (4.3) can be checked by an inductive argument similar to quad trees by using the expression

$$\pi_{n,j} = \frac{\binom{j}{t}\binom{n-1-j}{(m-1)(t+1)-1}}{\binom{n}{m(t+1)-1}}.$$

In particular, we can derive an $O$-transfer similar to Lemma 1 with the two numbers $\alpha, \rho$ there satisfying

$$\rho < \frac{(\alpha + t + 1) \cdots (\alpha + m(t+1) - 1)}{(t+2) \cdots (m(t+1))}.$$

Thus *the generalized $m$-ary search trees are also width-regular with the almost sure convergence* (4.5).

Note that $m$-ary search trees correspond to $t = 0$ and that $m = 2$ reduces to the so-called fringe-balanced or median-of-$(2t+1)$ binary search trees; see [8].

**7. Random increasing trees.** Bergeron, Flajolet and Salvy [2] proposed the increasing tree model: rooted trees with increasing labels along paths down from the root, such that the number of trees of a certain structure is prescribed in some general manner. The exponential generating function $\tau(z) := \sum_n \tau_{n \geq 1} z^n / n!$ for the number $\tau_n$ of increasing trees (of $n$ labels) has the form

$$(7.1) \qquad \tau'(z) = \phi(\tau(z)),$$

with $\tau(0) = 0$ and $\tau(1) = 1$ for some function $\phi(w)$ with $\phi(0) = 1$ and non-negative Taylor coefficients. The degree function $\phi(w)$ specifies how the trees are recursively formed. In this case, there are three representative varieties of increasing trees: (i) recursive trees with $\phi(w) = e^w$; (ii) binary increasing trees with $\phi(w) = (1 + w)^2$; (iii) plane-oriented recursive trees (PORTs) with $\phi(w) = 1/(1 - w)$.

We already studied the width of random recursive trees and random binary increasing trees (identically distributed as random binary search trees). We consider first PORTs and then mention other varieties of increasing trees (in some generality).

*Random PORTs.* PORTs are ordered (or plane) increasing trees without restriction on the degree of each node. To the best of our knowledge, such trees first appeared in a combinatorial form in [40], and in a more general probabilistic form in [44]; see also [2, 39, 43] for more details. Probabilistic properties of PORTs received much recent attention due partly to their close connection to random complex models; see [26] for more references.

The recurrence for the profile $Y_{n,k}$ is similar to (2.1), but with a very different underlying distribution (see [26])

$$Y_{n,k} \overset{d}{=} Y_{I_n,k-1} + Y^*_{n-I_n,k} \qquad (n \geq 2; k \geq 1),$$

where the $Y^*_{n,k}$'s are independent copies of $Y_{n,k}$ and

$$\pi_{n,j} = \mathbb{P}(I_n = j) = \frac{2\binom{2j-2}{j-1}\binom{2n-2j-2}{n-j-1}}{j\binom{2n-2}{n-1}} \qquad (1 \leq j < n).$$

We have

$$\Xi_n(u) = \frac{1}{(1+u)} \left( \frac{2\sqrt{\pi}}{\Gamma(u/2)} n^{(u+1)/2} + 1 \right)(1 + O(n^{-\varepsilon}))$$

uniformly for $|u| \leq C$ for any $C > 0$; see also [2, 31]. Then we have (4.1) with $f(u) = (u + 1)/2$, so that $f'(1) = \sigma^2 = 1/2$. Note that although the recurrence satisfied by $Y_{n,k}$ is not of the form (4.8), the technicalities are similar to those for recursive trees; see [26] for details. Thus *PORTS are width-regular with parameters $(1/2, 1/2)$; the almost sure convergence (4.5) also holds.*

The widths and profiles of random increasing trees for which $1/\phi(w)$ equals a polynomial also exhibit similar behaviors.

*Polynomial varieties.* We now show that the same width-regularity results (4.4) also hold for polynomial varieties of increasing trees; see [2]. Briefly, these are increasing trees in which each node has at most $d$ subtrees and the exponential generating function $\tau(z) := \sum_{n \geq 1} \tau_n z^n/n!$ satisfies (7.1) with

$$\phi(w) := \sum_{0 \leq j \leq d} \phi_j w^j \qquad (d \geq 2),$$

where $\phi_j \geq 0$ for $0 \leq j \leq d$ and $\phi_0, \phi_d > 0$. In this case, it is known that

$$
\begin{aligned}
(7.2) \quad \frac{\tau_n}{n!} &= \frac{p}{\Gamma(1/(d-1))}((d-1)\phi_d R)^{-1/(d-1)} \\
&\quad \times R^{-n} n^{-(d-2)/(d-2)} (1 + O(n^{-2/(d-1)})),
\end{aligned}
$$

where $p$ denotes the period of $\phi(v)$, $R := \int_0^\infty dv/\phi(v)$ and

$$\sum_{n,k} \mathbb{E}\{Y_{n,k}\} u^k \frac{z^n}{n!} = (\tau'(z))^{-u} \int_0^z (\tau'(v))^{1-u} \, dv;$$

see [2]. From these relations, we deduce the two estimates (4.1) and (4.2) with

$$
\begin{aligned}
f(u) &= \frac{d}{d-1} u - \frac{1}{d-1}, \\
g(u) &= \phi_d^{(1-u)/(d-1)} (R(d-1))^{(1-du)/(d-1)} \\
&\quad \times \frac{\Gamma(1/(d-1))}{\Gamma(du/(d-1))} \int_0^\infty \phi(v)^{-u} \, dv.
\end{aligned}
$$

Furthermore, the higher moments of $Y_{n,k}$ (centered or not) satisfy the recurrence

$$a_{n,k} = b_{n,k} + \sum_{1 \leq j < n} \varpi_{n,j} a_{j,k-1},$$

where

$$\varpi_{n,j} := \frac{(n-1)! \tau_j}{\tau_n j!} [z^{n-1-j}] \phi'(\tau(z)).$$

By (7.2), we then derive an $O$-transfer for $a_{n,k}$ similar to Lemma 1 with $\alpha$ and $\rho$ there satisfying

$$\rho < \frac{d-1}{dp^{d-1}}\left(\alpha + \frac{1}{d-1}\right);$$

from this the estimates (4.3) are then justified, implying the width-regularity properties (4.4).

*Random mobile trees.* These are increasing trees in which all subtrees are arranged in cyclic order and whose enumerating generating function satisfies (7.1) with $\phi(w) = 1 - \log(1 - w)$; see [2]. This example is less natural but very interesting because nodes are distributed in a rather different way and *the trees are not width-regular.*

First, the generating polynomial for the expected profile is given by

$$\Xi_n(u) = \sum_k \mu_{n,k} u^k = \frac{n!}{\tau_n} [z^n] \tau'(z)^u \int_0^z \tau'(v)^{1-u} \, dv.$$

Here the number $\tau_n$ of such trees satisfies

$$\frac{\tau_n}{n!} = R^{1-n} n^{-2} (1 + O(L_n^{-1})),$$

where $R = \int_0^\infty (1+v)^{-1} e^{-v} \, dv$. By singularity analysis (see [19]), we deduce that

$$\Xi_n(u) = g(u) n L_n^{u-1} \left( 1 + O\left( \frac{\log L_n}{L_n} \right) \right),$$

where the $O$-term holds uniformly for bounded complex $u$ and

$$g(u) = R^{-1} u \int_0^\infty e^{-v} (1+v)^{-u} \, dv.$$

Note that this is not of the form (4.1) and $g(1) = 1$. Thus such mobile trees are very "bushy" at each level (the root already having about $n/L_n$ nodes) and we have

$$\max_k \mu_{n,k} \sim \frac{n}{\sqrt{2\pi \log L_n}},$$

the mode being reached at $k \sim \log L_n$. The same methods of proof we used for recursive trees can be extended to show that

$$\mathbb{E}\{W_n\} \sim \frac{n}{\sqrt{2\pi \log L_n}},$$

a very different behavior from all types of random trees we have discussed.

**Acknowledgments.** This paper was prepared while both authors were visiting Institut für Stochastik und Mathematische Informatik, J. W. Goethe-Universität (Frankfurt). They thank Ralph Neininger and the Institute for their hospitality. We thank the referee for helpful suggestions (including, in particular, the term "width-regular") and comments.

SCHOOL OF COMPUTER SCIENCE
MCGILL UNIVERSITY
MONTREAL H3A 2K6
CANADA
E-MAIL: luc@cs.mcgill.ca

INSTITUTE OF STATISTICAL SCIENCE
ACADEMIA SINICA
TAIPEI 115
TAIWAN
E-MAIL: hkhwang@stat.sinica.edu.tw